# WHEN DO STEPWISE ALGORITHMS MEET SUBSET SELECTION CRITERIA?[1]

By Xiaoming Huo and Xuelei (Sherry) Ni

*Georgia Institute of Technology*

Recent results in *homotopy* and *solution paths* demonstrate that certain well-designed greedy algorithms, with a range of values of the algorithmic parameter, can provide solution paths to a sequence of convex optimization problems. On the other hand, in regression many existing criteria in subset selection (including $C_p$, AIC, BIC, MDL, RIC, etc.) involve optimizing an objective function that contains a counting measure. The two optimization problems are formulated as (P1) and (P0) in the present paper. The latter is generally combinatoric and has been proven to be NP-hard. We study the conditions under which the two optimization problems have common solutions. Hence, in these situations a stepwise algorithm can be used to solve the seemingly unsolvable problem. Our main result is motivated by recent work in *sparse representation*, while two others emerge from different angles: a direct analysis of sufficiency and necessity and a condition on the mostly correlated covariates. An extreme example connected with *least angle regression* is of independent interest.

**1. Introduction.** We consider two types of optimization problem:

- an optimization problem that is based on a counting measure,

$$\text{(P0)} \qquad \min_x \|y - \Phi x\|_2^2 + \lambda_0 \cdot \|x\|_0,$$

where $\Phi \in \mathbb{R}^{n \times m}, x \in \mathbb{R}^m, y \in \mathbb{R}^n$, the notation $\|\cdot\|_2^2$ denotes the sum of squares of the entries of a vector, the constant $\lambda_0 \geq 0$ is an algorithmic parameter and the quantity $\|x\|_0$ is the number of nonzero entries in the vector $x$;

- an optimization problem that depends on a sum of absolute values,

$$\text{(P1)} \qquad \min_x \|y - \Phi x\|_2^2 + \lambda_1 \cdot \|x\|_1,$$

Received July 2005; revised July 2006.

[1]Supported by in part NSF Grant DMS-01-40587.

*AMS 2000 subject classification.* 62J07.

*Key words and phrases.* Subset selection, model selection, stepwise algorithms, convex optimization, concurrent optimal subset.

where $\|x\|_1 = \sum_{i=1}^{m} |x_i|$ for the vector $x = (x_1, x_2, \ldots, x_m)^T$ and where the constant $\lambda_1 \geq 0$ is another algorithmic parameter.

Note that $\|x\|_0$ (resp., $\|x\|_1$) is a quasi-norm (resp., norm) in $\mathbb{R}^m$. In the literature on *sparse presentation*, these are called the $\ell_0$-norm and the $\ell_1$-norm, respectively. The notation (P0) and (P1) also appears in [8], with slightly different definitions.

In subset selection under linear regression, many well-known criteria—including the $C_p$ statistic, the Akaike information criterion (AIC), the Bayesian information criterion (BIC), minimum description length (MDL), the risk inflation criterion (RIC) and so on—are special cases of (P0), resulting from the assignment of different values to $\lambda_0$. It is shown in this paper that problem (P0) is, in general, NP-hard (Theorem 2.1). The NP-hardness has been known for many years, but to the best of our knowledge, no paper has formally presented a proof of this yet.

At the same time, (P1), which has a long history that will be reviewed later, is the mathematical problem that is called upon in [45]. Recent advances (details and references are provided in Section 2.2) demonstrate that some stepwise algorithms (e.g., [10, 38, 39]) reveal the solution paths of problem (P1) while the parameter $\lambda_1$ takes a range of values. More importantly, most of these algorithms take only a polynomial number of operations (i.e., they are polynomial-time algorithms). In fact, the complexity of finding a solution path for (P1) is the same as that of implementing an ordinary least square fit [10].

The major objective of this paper is to find out when (P0) and (P1) give the same result in subset selection. A subset that corresponds to the nonzero subset of the minimizer of (P0) [resp., (P1)] is called a *type-0* (resp., *type-1*) *optimal subset with respect to* $\lambda_0$ (resp., $\lambda_1$). A subset that is both type-0 and type-1 optimal is called a *concurrent optimal subset*. It is known that there is a necessary and sufficient condition for the type-1 optimal subset and that this condition can be verified in polynomial time. However, there is no polynomial-time necessary-and-sufficient condition for the type-0 optimal subset. We search for easily verified (i.e., polynomial-time) sufficient conditions for type-0 optimal subsets. When sufficient conditions are available, given solutions of (P1) by a stepwise algorithm, we can determine whether (P0) has been solved. The title of this paper reflects such an objective.

The main contributions of this paper are two verifiable sufficient conditions for (P0), Theorems 3.1 and 3.2. The latter is an improved version of the former. Other conditions are generally standard and known. They are presented subsequently for the sake of completeness.

The paper is organized as follows. Section 2 reviews the subset selection criteria that can be formulated as (P0), as well as the literature on (P1). Two cases are studied/reviewed. Section 3 contains the main results. Section [4](4) presents some associated conditions, which are either known or relatively easy to verify. Section [5](5) discusses some related work. A brief conclusion is provided in Section [6](6). Proofs are relegated to the appendices when convenient.

**2. Formulation and review of literature.** We consider a regression setting. Let $\Phi \in \mathbb{R}^{n \times m}$ ($n > m$) denote a model matrix. Vectors $x \in \mathbb{R}^m$ and $y \in \mathbb{R}^n$ are coefficient and response vectors. The columns of matrix $\Phi$ are *covariates*. A regression model is $y = \Phi x + \varepsilon$, where $\varepsilon$ is a random vector. Let $\mathbf{I} = \{1, 2, \ldots, m\}$ denote the set of indices of the coefficients. A subset of coefficients (or covariates) is denoted by $\Omega$ ($\Omega \subseteq \mathbf{I}$). Let $|\Omega|$ denote the cardinality of the set $\Omega$. Let $x_\Omega$ denote the coefficient vector that takes only nonzero values when the coefficient indices are in the subset $\Omega$. To choose a subset $\Omega$, a subset selection problem has two competing objectives: (1) the residuals, $y - \Phi x_\Omega$, are close to zero and (2) the size of the set $\Omega$ is small.

2.1. *Subset selection criteria and* (P0). There exists an extensive body of literature on the criteria regarding subset selection. Miller [31], Burnham and Anderson [2] and George [19] all give excellent reviews. An interesting fact is that a majority of these criterion can be unified under (P0), where $\|y - \Phi x\|_2^2$ is the residual sum of squares [denoted by RSS($x$)] and where the constant $\lambda_0$ depends on the criterion. Some well-known results are the Akaike information criterion (AIC) [1], $C_p$ [20, 30], the Bayesian information criterion (BIC) [44], minimum description length (MDL) (see the equivalence between BIC and MDL in [25], Section 7.8), the risk inflation criterion (RIC) [15] and so on. We refer to George [19] for the details. In this paper, the "subset selection criteria" that appears in the title encompasses all of the foregoing criteria.

Solving (P0) generally requires an exhaustive search of all subsets. When $\|x\|_0$ (i.e., the number of covariates) increases, the methods based on exhaustive search rapidly become impractical. Innovative ideas have been developed to reduce the number of subsets being searched; see [17, 32], as well as some later improvements, [18, 35, 36, 40, 41]. All of these methods adopt a branch-and-bound (B&B) strategy. Improvements can be achieved by modifying the structure in B&B or by applying stronger optimality tests. Despite these efforts, when the number of covariates ($m$) is moderately large (e.g., $m = 50$), the subset search cannot generally be carried out, unless the model matrix $\Phi$ possesses some special structure.

In fact, solving (P0) is an NP-hard problem! The following theorem can be considered as an extension of a result originally presented in [33]. The proof of the theorem appears in Appendix A.1.

THEOREM 2.1. *Solving* (P0) *with a fixed $\lambda_0$ is an NP-hard problem.*

2.2. *Stepwise algorithms and* (P1). Due to the difficulty of solving (P0), a *relaxation* idea has been proposed. The relaxation replaces the $\ell_0$-norm with the $\ell_1$-norm in the objective, which leads to (P1).

Santosa and Symes [43] is considered the first modern appearance of the formulation (P1). The idea of relaxation has been studied extensively in the literature on *sparse representation*. Some representative papers are (roughly in chronological order) [4, 6, 7, 8, 11, 16, 22, 23, 29, 46, 47], and so on. A full review is well beyond the scope of this paper. The problem of sparse representation has a different emphasis, involving the derivation of a priori conditions instead of a posteriori conditions, as in the present paper.

At the same time, (P1) has been proposed in the statistics literature as a method of subset selection. This has been termed the *Lasso method* [45]. An interesting recent development—least angle regression (LARS) [10]— demonstrates that certain stepwise algorithms can reveal the solutions to (P1) with varying values of $\lambda_1$, based on the idea of *homotopy* (see [38]). More recent analysis further demonstrates that stepwise algorithms can literally render the entire solution path in a large class of problems; see [24] and the references therein. The homotopy continuation method [39] and the *subdifferential* are the key technical tools in this development. [42] and [37] are useful references.

2.3. *Case studies.* We present two cases that have been instructive to us.

2.3.1. *An extreme example.* We construct an extreme example, in which a sophisticated stepwise algorithm (LARS) can miss an optimal subset. This example has played an inspirational role in our study of the equivalence conditions, which are discussed further in the next section. Details of the LARS algorithm can be found in [10], Section 2. We believe that this example is interesting in its own right.

The example is constructed as follows. Let $\phi_i \in \mathbb{R}^n, i = 1, 2, \ldots, m$, denote the $i$th column of the model matrix $\Phi$. Hence, $\Phi = [\phi_1, \phi_2, \ldots, \phi_m]$. Let $\delta_i \in \mathbb{R}^n, i = 1, 2, \ldots, m$, denote the Dirac vector taking 1 at the $i$th position and zero elsewhere. For $i = m - A + 1, m - A + 2, \ldots, m$, let $\phi_i = \delta_i$, where $A$ is a positive integer. Consider a signal $s = \frac{1}{\sqrt{A}} \sum_{i=m-A+1}^{m} \phi_i$. Obviously, in this case the optimal subset is $\{m - A + 1, \ldots, m\}$. For the first $m - A$ columns of $\Phi$, let $\phi_j = a_j \cdot s + b_j \cdot \delta_j$, where $1 \leq j \leq m - A$ and $a_j^2 + b_j^2 = 1$. Note that the $\phi_i$'s and $s$ are all unit-norm vectors. Hereafter, for simplicity, we always assume that $1 \leq j \leq m - A$ and $m - A + 1 \leq i \leq m$. It is easy to verify that

$$\langle s, \phi_j \rangle = a_j \quad \text{and} \quad \langle s, \phi_i \rangle = 1/\sqrt{A}.$$

In this example, we choose $1 > a_1 > a_2 > \cdots > a_{m-A} > 1/\sqrt{A} > 0$.

THEOREM 2.2 (An extreme example). *In the above example, the LARS algorithm chooses covariates $\phi_1, \phi_2, \ldots, \phi_{m-A}$ one-at-a-time and by the same order in the first $m - A$ steps.*

It takes some effort to verify the above theorem. We refer to the technical report [26], which is a longer version of this paper, as well as to the thesis [34].

Readers may notice that in the above example, the covariates are *not* standardized, while in the LARS algorithm, choosing covariates according to the inner product implies the standardization of covariates. A discussion in [26], Theorem 3.4, shows that this can be remedied by an orthogonal transformation.

The foregoing example is developed in a fairly general form, with controlling parameters $A$ and $m$. To illustrate how dramatic this example can be, let us consider the case where $A = 10$ and $m = 1{,}000{,}000$. Based on the previous description, the LARS algorithm will select the first 999,990 covariates before it selects any of the last ten covariates. At the same time, the optimal subset is formed by the last ten covariates. Another example regarding the performance of LARS can be found in [48], which has a different emphasis.

This example is motivated by an early example in [4], which can be traced further back to [3] and [5] in the analysis of some stepwise algorithms (e.g., orthogonal matching pursuit) in signal processing. Our example is similar in spirit; however it is different in constructional details.

2.3.2. *Subset selection with orthogonal model matrix.* The following result is well known: for an orthogonal model matrix $\Phi$, when $\sqrt{\lambda_0} = \lambda_1/2$, solutions to (P0) and (P1) have the same support. Moreover, at each position the solutions differ by a constant $\lambda_1/2$. A partial list of references for such a result includes [10, 38, 45], and many more. For readers who are familiar with soft-thresholding and hard-thresholding [9], this result should not come as a surprise.

The above two examples collectively motivate us to pursue sufficient conditions that guarantee common support in the solutions of (P0) and (P1).

**3. Main results.** A general sufficient condition for (P0) is derived. It is motivated by a recent approach which has appeared in applied mathematics; see [22]. We have modified their approach to solve a different mathematical problem.

Recall that $x \in \mathbb{R}^m$ denotes a coefficient vector. Denote the corresponding residual vector by $\varepsilon = y - \Phi x$. Recall that $y \in \mathbb{R}^n$ and $\Phi \in \mathbb{R}^{n \times m}$ are the response vector and the model matrix, respectively. Let $\Omega$ denote the support of the vector $x : \Omega = \text{supp}(x)$. For an integer $k \geq 1$, let

$$\sigma_{\min,k}^2 = \inf \frac{\|\Phi \delta\|_2^2}{\|\delta\|_2^2} \quad \text{subject to} \quad \|\delta\|_0 \leq k.$$

The above quantity reflects a certain property of the model matrix. Furthermore, for a vector $v \in \mathbb{R}^n$ and an integer $k \geq 1$, we define

$$c(v, k) = \sqrt{\sum_{i=1}^{k} v_{(i)}^2},$$

where $|v_{(1)}| \geq |v_{(2)}| \geq \cdots \geq |v_{(n)}|$ are the nonincreasing-ordered magnitudes of the entries of the vector $v$. For finite $k$, we assume that the quantities $c^2(\Phi^T \varepsilon, k)$ and $\sigma_{\min,k}^2$ have been computed. The following theorem provides a sufficient condition for a subset to be included in a type-0 optimal subset with respect to $\lambda_0$.

THEOREM 3.1 (Main result 1). *A subset of coefficients $\Omega$ is given. Suppose that coefficient vector $x$ is the minimizer of the function $\|y - \Phi x\|_2^2$ subject to $\mathrm{supp}(x) \subset \Omega$. Let $\varepsilon = y - \Phi x$.*

*(1) If $\min_{i \in \Omega} |x_i| > q_1(|\Omega|)$, then, with respect to $\lambda_0$, there is no type-0 optimal subset whose support is of size less than $|\Omega|$.*

*(2) Furthermore, if $\min_{i \in \Omega} |x_i| > q(|\Omega|)$, then, with respect to $\lambda_0$, we have $\Omega \subset \Omega'$, where $\Omega'$ is the type-0 optimal subset with respect to $\lambda_0$.*

*The quantities $q_1(\cdot)$ and $q(\cdot)$ are defined as follows. For an integer $k \geq 1$,*

$$q_1(k) = \sup_{m < k} \frac{c(\Phi^T \varepsilon, 1) + \sqrt{c^2(\Phi^T \varepsilon, k+m) + (k-m)\lambda_0 \sigma_{\min,k+m}^2}}{\sigma_{\min,k+m}^2},$$

$$q_2(k) = \sup_{m \geq k} \frac{c(\Phi^T \varepsilon, 1) + \sqrt{c^2(\Phi^T \varepsilon, k+m) + (k-m)\lambda_0 \sigma_{\min,k+m}^2}}{\sigma_{\min,k+m}^2}$$

*and*

$$q(k) = \max\{q_1(k), q_2(k)\}.$$

The proof of this theorem appears in Appendix A.2.

Note that quantities $q_1(\cdot)$ and $q_2(\cdot)$ have the same objective function. However, the ranges of the variable $m$ are different. Because $q_1(k)$ requires only a finite choice of the variable $m$ (recall that $m < k$), it is computable. It is not straightforward to show that for any $k \geq 1$, the quantity $q_2(k)$ will exist. In this paper, we assume the existence of this quantity.

Readers may compare the above with the test proposed in [22]. That test is related to the optimality in sparse representations.

In Theorem 3.1, the quantities $q_1(\cdot)$ and $q(\cdot)$ require multiple values of $\sigma_{\min,k}^2$, for a range of values of $k$. Compared to the quantities $c(\cdot, k)$, it is harder to compute the $\sigma_{\min,k}^2$'s. Inspired by the derivation in Theorem 2 of [22], we derive a sufficient condition, which depends only on $\sigma^2_{\min,|\Omega|}$, where $\Omega$ is the subset being tested. To state our result, the following quantity needs to be defined: for an integer $m \geq 1$ and a given integral constant $M$, let

$$\lambda(m; M) = 1 - \frac{M}{\sqrt{m}} \sup_{|\mathcal{I}| \leq m} \sup_{k \notin \mathcal{I}} \|\Phi_{\mathcal{I}}^+ \phi_k\|_2,$$

where $\mathcal{I}$ is a subset of indices, $|\mathcal{I}|$ denotes the size of this subset, the matrix $\Phi_{\mathcal{I}}$ is a submatrix of $\Phi$ whose column indices form the set $\mathcal{I}$, $\Phi_{\mathcal{I}}^+ = (\Phi_{\mathcal{I}}^* \Phi_{\mathcal{I}})^{-1} \Phi_{\mathcal{I}}^*$ is the Moore–Penrose pseudo-inverse [21], with $(\cdot)^*$ denoting the adjoint, and $\phi_k$ is the $k$th column (i.e., covariate) in $\Phi$. Given $m$, the quantity $\lambda(m)$ can be computed by enumerating all $m$-subsets of the covariates.

Now we present another sufficient condition.

THEOREM 3.2 (Main result 2). *A subset of coefficients $\Omega$ is given. Suppose that coefficients vector $x$ is the minimizer of the function $\|y - \Phi x\|_2^2$ subject to $\mathrm{supp}(x) \subset \Omega$. Suppose it is known a priori that the size of the type-0 optimal subset is no larger than $M$. If $\min_i |x_i| > q'(|\Omega|, M)$, then the set $\Omega$ is at least a subset of the type-0 optimal subset. Here, the quantity $q'(\cdot)$ is defined, for integer $k \geq 1$ and constant $M$, as*

$$q'(k, M)$$
$$= \sup_{1 \leq m \leq M} \left( c(\Phi^T \varepsilon, 1) + \sqrt{c^2(\Phi^T \varepsilon, k) + \lambda_0 \cdot \frac{k^2(k-m)}{(k+m)^2} \cdot \sigma^2_{\min,k} \cdot \lambda^2(k; m)} \right)$$
$$\times \left( \frac{k}{(k+m)} \sigma^2_{\min,k} \cdot \lambda^2(k; m) \right)^{-1}.$$

See the proof in Appendix A.3.

If the model matrix $\Phi$ is orthonormal, readers can verify that $\sigma^2_{\min,k} = 1$ and $\lambda(m; M) = 1$. This brings about significantly simplified criteria in Theorem 3.1 and Theorem 3.2. Compared with the case when the model matrix is orthogonal, the new criteria are less attractive. We consider this a price to be paid for the generality.

The two results here focus on the type-0 optimal subset. Given a type-1 optimal subset (which can be derived from some efficient algorithm), one can easily calculate the least square estimator according to it and use this estimator and subset to test whether the subset is also type-0 optimal.

**4. Other conditions of equivalence.** In Section 4.1, we give a sufficient and necessary condition for a subset to be the concurrent optimal subset. Checking this condition cannot be achieved in polynomial time [recall that (P0) is NP-hard]. In Section 4.2, we ask when the $k$ most correlated covariates form the concurrent optimal subset. A sufficient condition is derived. This result is easy to check, but too restrictive.

4.1. *Sufficient and necessary conditions.* Before moving on to the specific discussion, we introduce a sufficient and necessary condition for a concurrent optimal subset. Let $I_1$ denote a subset of indices. Let $\Phi_1$ and $x_1$ denote columns of $\Phi$ and entries of $x$ with indices from $I_1$. Let $\Phi = [\Phi_1 \ \Phi_2]$. Here, a permutation that does not change the problem is implied. The following can easily be verified.

LEMMA 4.1 [Sufficient and necessary condition for (P0)]. $I_1$ *is the optimal subset of* (P0) *if and only if the value*

$$(1) \qquad y^T y - y^T \Phi_1 (\Phi_1^T \Phi_1)^{-1} \Phi_1^T y + \lambda_0 \cdot \|x_1\|_0$$

*is the minimum of the objective in* (P0).

The following is well known (see, e.g., [29, 38, 47]).

LEMMA 4.2 [Sufficient and necessary condition for (P1)]. $I_1$ *is the optimal subset of* (P1) *if and only if there exists a vector $\omega$ such that*

$$(2) \qquad \Phi^T y = \begin{pmatrix} \Phi_1^T \Phi_1 \\ \Phi_2^T \Phi_1 \end{pmatrix} x_1 + \begin{pmatrix} \frac{\lambda_1}{2} \cdot \text{sign}(x_1) \\ \omega \end{pmatrix}$$

*holds and* $\|\omega\|_\infty \le \lambda_1/2$.

The following can be easily derived from the above two lemmas.

COROLLARY 4.3 [Sufficient and necessary condition (for concurrence)]. $I_1$ *is the concurrent optimal subset of* (P0) *and* (P1) *if and only if* (1) *and* (2) *are true. Moreover, with $\widetilde{x}_0$ and $\widetilde{x}_1$ the solutions of* (P0) *and* (P1), *respectively, we have*

$$(3) \qquad (\widetilde{x}_0 - \widetilde{x}_1)_{I_1} = (\Phi_1^T \Phi_1)^{-1} \cdot \frac{\lambda_1}{2} \cdot \text{sign}((\widetilde{x}_1)_{I_1}).$$

For the equation above, consider

$$(\widetilde{x}_0)_{I_1} = (\Phi_1^T \Phi_1)^{-1} \Phi_1^T y$$

and

$$\Phi_1^T y = (\Phi_1^T \Phi_1)(\widetilde{x}_1)_{I_1} + \frac{\lambda_1}{2} \cdot \text{sign}((\widetilde{x}_1)_{I_1}).$$

By combining the above two, equation (3) follows.

The above theorem gives a necessary and sufficient condition for a concurrent optimal subset. Further comments follow.

REMARK 4.4. Equation (3) provides a method for computing $\widetilde{x}_1$, given that $\widetilde{x}_0$ is available and represents the optimal solution. Evidently,

$$(\widetilde{x}_1)_{I_1} = (\widetilde{x}_0)_{I_1} - \frac{\lambda_1}{2}(\Phi_1^T \Phi_1)^{-1} \cdot \text{sign}((\widetilde{x}_1)_{I_1}).$$

REMARK 4.5. Note that

$$\begin{aligned}\Phi(\widetilde{x}_0 - \widetilde{x}_1) &= \Phi_1(\widetilde{x}_0 - \widetilde{x}_1)_{I_1} \\ &= \frac{\lambda_1}{2} \cdot \Phi_1(\Phi_1^T \Phi_1)^{-1} \cdot \text{sign}((\widetilde{x}_1)_{I_1}),\end{aligned}$$

which is an equiangular vector among the columns of $\Phi_1$. Hence, when optimality is achieved in both (2) and (3), the difference between the two predicted vectors is an equiangular vector.

4.2. *A sufficient condition for mostly correlated covariates.* We introduce a set of sufficient conditions which depend only on the correlations between the response $y$ and the covariates $\phi_i$, as well as the maximum correlation between the covariates. For simplicity, we assume that the response $y$ and the covariates $\phi_i$ are all standardized. It is not hard to see that $|\langle y, \phi_i \rangle| \leq 1$, $i = 1, 2, \ldots, m$, and $|\langle \phi_i, \phi_j \rangle| \leq 1, 1 \leq i, j \leq m$. Denote $z = \Phi^T y = (z_1, z_2, \ldots, z_m)^T$. Without loss of generality, we assume $|z_1| > |z_2| > \cdots > |z_m|$. We want to find sufficient conditions such that the subset $A_1 = \{\phi_1, \phi_2, \ldots, \phi_k\}$ is the solution to both (P0) and (P1): the $k$ most correlated covariates (with the response) form the optimal subset. Clearly, an optimal subset does not need to consist of the most correlated covariates with the response. Due to this additional condition, this set of conditions is *restrictive*. The restrictiveness is illustrated in an example in Section 4.2.1.

Denote

$$\mu = \max_{\substack{1 \leq i,j \leq m \\ i \neq j}} |\langle \phi_i, \phi_j \rangle|.$$

We have the following theorem.

THEOREM 4.6. *For a given $\lambda_0$ and correlations $z_1, z_2, \ldots, z_k$, if the three conditions*

$$(4) \qquad [1 - (k-1)\mu]z_k^2 \geq 2(k-1)^2\mu + z_{k+1}^2[1 + (k-1)\mu],$$

$$(5) \qquad z_{k+1}^2 \leq \lambda_0(1 - \Delta) - \frac{(2k-1)\mu}{1 + (k-1)\mu}\sum_{i=1}^k z_i^2,$$

$$(6) \qquad z_k^2 \geq \lambda_0 + \frac{(2k-3)\mu}{1 + (k-1)\mu}\sum_{i=1}^k z_i^2$$

*are satisfied, where* $\Delta = n \cdot \mu$ *in* (5), *then the subset* $A_1$ *is the type-*0 *optimal subset.*

To prove the above theorem, we can show that for subsets of size equal to $k$, greater than $k$, less than $k$, the above three conditions guarantee that subset $A_1$ is the type-0 optimal subset. A detailed proof can be found in [26] or [34]. Anyone whose interest is restricted to (P0) should now be satisfied. The following is to establish a condition for concurrent optimality.

REMARK 4.7. Conditions (4), (5) and (6) are independent, that is, none of them can be derived from the other two.

The following theorem states the condition for the set $A_1 = \{\phi_1, \phi_2, \ldots, \phi_k\}$ to be the type-1 optimal subset; see the proof in [26] or [34].

THEOREM 4.8. *Given* $\lambda$ *and* $k$, *if*

$$(7) \qquad \frac{\lambda}{2} - |z_{k+1}| \geq \frac{\sqrt{k}\mu}{1-(k-1)\mu}\sqrt{\sum_{i=1}^{k}\left(|z_i| + \frac{\lambda}{2}\right)^2},$$

*then subset* $A_1$ *is the type-*1 *optimal subset.*

The following corollary gives a sufficient condition for $A_1$ to be the concurrent optimal subset.

COROLLARY 4.9. *Given conditions* (4), (5), (6) *and* (7), *subset* $A_1$ *is the concurrent optimal subset.*

4.2.1. *Restrictiveness of the aforementioned sufficient conditions.* Readers may have noticed that the four conditions in the previous section are restrictive. One can easily find an example that does not satisfy these conditions, but which still has the concurrent optimal subset $A_1$.

An example can be established as follows. Suppose that $n, m$ and $k$ are three positive integers satisfying $n > m > k$ and $n \geq m + k$. Let $a_i$ denote the $i$th entry of vector $\mathbf{a} \in \mathbb{R}^k$ with $|a_1| \geq |a_2| \geq \cdots \geq |a_k|$. Let $I_{m \times m} \in \mathbb{R}^{m \times m}$ be an identity matrix and $\Phi_a \in \mathbb{R}^{k \times k}$ be the diagonal matrix with the $i$th diagonal entry being equal to $a_i$. Consider

$$\Phi = \text{standardized}\left\{\begin{pmatrix} \Phi_a\ \mathbf{0}_{k\times(m-k)} \\ I_{m\times m} \\ \mathbf{0}_{(n-k-m)\times m} \end{pmatrix}\right\}, \qquad y = \sum_{i=1}^{k}\phi_i,$$

where "standardized$\{M\}$" refers to the standardization of all of the columns of the matrix $M$, the matrices $\mathbf{0}_{k \times (m-k)}$ and $\mathbf{0}_{(n-k-m) \times m}$ consist entirely of zeros and $\phi_i$ is the $i$th column of $\Phi$. The optimal solution consists of the first $k$ covariates and these covariates have larger correlations with $y$. However, there are many choices of $m, n, k$ and the vector $\mathbf{a}$ with which condition (4) is not satisfied. As a special case, consider the following simple example: $n = 10, m = 7, k = 3$ and $\mathbf{a} = (-1\ 1\ 0)^T$. It is not hard to verify that $\mu(\Phi) = 0.1667, z_3 = 0.7379, z_4 = -0.3162, [1 - (k-1)\mu]z_k^2 = 0.3630$ and $2(k-1)^2\mu + z_{k+1}^2[1 + (k-1)\mu] = 0.9117$. Hence, (4) does not hold in this case.

**5. Discussion.** The question addressed in this paper has a unique aspect. We have the following application in mind: supposing a stepwise algorithm finds a path of type-1 optimal subsets, then given verifiable (polynomial-time) conditions that are derived in this paper, one knows whether a type-0 optimal subset has been found. As mentioned earlier, our results potentially facilitate polynomial-time solutions to seemingly NP-hard problems.

Our problem is different from that of analyzing statistical properties of the estimators. These properties include consistency, rate of convergence, asymptotic normality and so on. We found the oracle properties derived in [12] very interesting. However, Fan and Li [12] do not address whether their estimator—smoothly clipped absolute deviation penalty (SCAD)—can be computed in polynomial time. In fact, because of the possible exponential number of local optima, it is strongly believed that SCAD cannot be solved in polynomial time. Hence, an interesting question will be: when can one verify that SCAD is indeed solved by a polynomial-time algorithm? That is, we want to derive some sufficient conditions similar to those in the present paper. Note that Fan and Peng [14] give a fundamental description of when oracle properties (as well as other properties) are achievable, while a recent manuscript by Zou [50] proves the oracle property for a method that is rooted in the Lasso.

As pointed out by an anonymous referee, there are two categories of equivalence conditions for (P0) and (P1): *a priori* conditions determine in advance when solving (P1) will identify a solution to (P0), while *a posteriori* conditions take a given subset of covariates (produced in any manner) and determine whether it is an optimal subset for (P0). The main results in this paper belong to the latter class. Given the target application described at the beginning of this section, it is not surprising that the latter is more interesting to us than the former. Moreover, a subset satisfying the former will most likely satisfy the latter, which implies that the a posteriori conditions are more powerful in the target application because they can identify more cases of equivalence.

Subset selection has applications in feature selection. There are two major approaches in feature selection: *filter* and *wrapper*; see [27, 28, 32] for details. Our formulations are closely related to wrappers. A recent survey paper by Fan and Li [13] gives an excellent overview of the statistical challenges associated with high-dimensional data, including feature selection and feature extraction. Besides many contemporary applications, as summarized in [13], other applications are foreseeable. For example, subset selection is a critical problem in supersaturated design. A citation search of Wu [49] will provide most of the existing literature. A numerically efficient condition on the optimality of subsets has the potential to identify a good design.

**6. Conclusion.** Stepwise algorithms can be numerically efficient, that is, polynomial-time. Specially designed stepwise algorithms can find type-1 optimal subsets in subset selection. We have derived sufficient conditions to test whether these type-1 optimal subsets are also type-0 optimal. Such an approach allows polynomial-time algorithms to locate concurrent optimal subsets, which, otherwise, generally requires solving an NP-hard optimization problem.

## APPENDIX A: PROOFS

### A.1. Proof of Theorem 2.1. Let

$$f(m) = \min_{x\,:\,\|x\|_0 \leq m} \|y - \Phi x\|_2^2,$$

where all of the symbols are defined in (P0). It is evident that the point array $(m, f(m))$, $m = 1, 2, \ldots$, forms a nonincreasing curve in the positive quadrant.

We first establish the existence of an integer $m_0$, such that value $f(m_0) + \lambda_0 m_0$ minimizes the objective in (P0). Note that there are a finite number of $m$'s such that $\lambda_0 m \leq f(1) + \lambda_0 \cdot 1$. This inequality gives an upper bound on $m$'s that satisfy $f(m) + \lambda_0 m \leq f(1) + \lambda_0 \cdot 1$. Among this finite number of $m$'s, there is at least one $m_0$ that minimizes the value of the function $f(m) + \lambda_0 m$.

Define $\varepsilon = f(m_0)$. In general, we can assume that $\varepsilon > 0$ because if $\varepsilon = 0$, then the response $y$ can be superposed by a small (more specifically, no more than $m_0$) number of columns of the matrix $\Phi$, which is a special case.

Using the idea of the Lagrange multiplier, we can see that solving (P0) with $\lambda_0$ is equivalent to solving the sparse approximate solution (SAS) problem in [33], Section 2, with $\varepsilon$, which is proved in [33] to be NP-hard. Hence, in general, solving (P0) is NP-hard.

### A.2. Proof of Theorem 3.1. Suppose that $\Omega'$ is the type-0 optimal subset, with corresponding coefficient vector $x'$. We must have

(8) $$\|y - \Phi x'\|_2^2 + \lambda_0 \|x'\|_0 \leq \|y - \Phi x\|_2^2 + \lambda_0 \|x\|_0.$$

Denoting $\delta = x' - x$, we have $\|\delta\|_0 \leq |\Omega| + |\Omega'|$. We will prove that

(9) $\qquad\qquad$ if $|\Omega'| < |\Omega|$, $\qquad$ then $\|\delta\|_\infty \leq q_1(\Omega)$

and

(10) $\qquad\qquad$ for any $\Omega'$, $\qquad \|\delta\|_\infty \leq q(\Omega)$.

To see the above, a reformulation of (8) gives

$$\|\varepsilon - \Phi\delta\|_2^2 \leq \|\varepsilon\|_2^2 + \lambda_0(|\Omega| - |\Omega'|),$$

which is equivalent to

(11) $\qquad\qquad \|\Phi\delta\|_2^2 \leq 2\langle \Phi^T \varepsilon, \delta\rangle + \lambda_0(|\Omega| - |\Omega'|),$

where $\langle \cdot, \cdot \rangle$ denotes the inner product between two vectors. Define $\delta' = \sigma_{\min,|\Omega|+|\Omega'|}^2 \cdot \delta$. Because $\|\Phi\delta\|_2^2 \geq \sigma_{\min,|\Omega|+|\Omega'|}^2 \|\delta\|_2^2$ and (11) hold, we have

$$\|\delta'\|_2^2 \leq 2\langle \Phi^T\varepsilon, \delta'\rangle + \lambda_0(|\Omega| - |\Omega'|) \cdot \sigma_{\min,|\Omega|+|\Omega'|}^2.$$

The above is equivalent to

$$\|\Phi^T\varepsilon - \delta'\|_2^2 \leq \|\Phi^T\varepsilon\|_2^2 + \lambda_0(|\Omega| - |\Omega'|) \cdot \sigma_{\min,|\Omega|+|\Omega'|}^2.$$

Define $\varepsilon^* = \Phi^T\varepsilon$. The above inequality leads to

$$\sum_{i\in\Omega\cup\Omega'}(\varepsilon_i^* - \delta_i')^2 \leq \sum_{i\in\Omega\cup\Omega'}(\varepsilon_i^*)^2 + \lambda_0(|\Omega| - |\Omega'|) \cdot \sigma_{\min,|\Omega|+|\Omega'|}^2.$$

The above immediately leads to

$$\sup_{i\in\Omega\cup\Omega'} |\delta_i'| \leq \sup_{i\in\Omega\cup\Omega'} |\varepsilon_i^*| + \sqrt{\sum_{i\in\Omega\cup\Omega'}(\varepsilon_i^*)^2 + \lambda_0(|\Omega| - |\Omega'|) \cdot \sigma_{\min,|\Omega|+|\Omega'|}^2}.$$

Dividing both sides by $\sigma_{\min,|\Omega|+|\Omega'|}^2$, we have

(12)
$$\sup_{i\in\Omega\cup\Omega'} |\delta_i|$$
$$\leq \frac{c(\Phi^T\varepsilon, 1) + \sqrt{c^2(\Phi^T\varepsilon, |\Omega| + |\Omega'|) + \lambda_0(|\Omega| - |\Omega'|) \cdot \sigma_{\min,|\Omega|+|\Omega'|}^2}}{\sigma_{\min,|\Omega|+|\Omega'|}^2}.$$

Recalling the definitions of $q_1(\cdot)$ and $q(\cdot)$, (9) and (10) can be derived directly from (12).

We are now able to verify item (1) of the theorem. Suppose that there is a type-0 optimal subset $\Omega'$ satisfying $|\Omega'| < |\Omega|$. We have

$$|x_i'| \geq |x_i| - |x_i - x_i'| \geq |x_i| - q_1(\Omega) > 0.$$

The second inequality is based on (9) and the last inequality follows from the condition in item (1). The above implies $\Omega \subset \Omega'$, which contradicts $|\Omega'| < |\Omega|$. We have proven item (1).

The proof of item (2) is very similar to the proof of (1). We omit the obvious details.

**A.3. Proof of Theorem 3.2.** The beginning of the proof is the same as that of the proof of the previous theorem. It begins to deviate at stage (11). For the reader's convenience, we restate inequality (11):

$$\|\Phi\delta\|_2^2 \leq 2\langle \Phi^T\varepsilon, \delta\rangle + \lambda_0(|\Omega| - |\Omega'|). \tag{13}$$

Readers are referred to the previous proof for the meanings of the notation.

First, we have

$$\langle \Phi^T\varepsilon, \delta\rangle \leq \sum_{i=1}^{n} |b_{(i)}| \cdot |\delta_{(i)}|, \tag{14}$$

where $|\delta_{(1)}| \geq |\delta_{(2)}| \geq \cdots \geq |\delta_{(n)}|$ is the ordered list of the magnitudes of the entries in the vector $\delta$. Similarly, $|b_{(1)}| \geq |b_{(2)}| \geq \cdots \geq |b_{(n)}|$ is the ordered list of the magnitudes of the entries in the vector $\Phi^T\varepsilon$. We denote $\Phi^T\varepsilon$ by $b$. The following manipulations are needed:

$$\begin{aligned}
\text{R.H.S. of (14)} &= \sum_{i=1}^{|\Omega|} |b_{(i)}| \cdot |\delta_{(i)}| + \sum_{i=|\Omega|+1}^{n} |b_{(i)}| \cdot |\delta_{(i)}| \\
&\leq \sum_{i=1}^{|\Omega|} |b_{(i)}| \cdot |\delta_{(i)}| + |b_{(|\Omega|+1)}| \cdot \sum_{i=|\Omega|+1}^{n} |\delta_{(i)}| \\
&\leq \sum_{i=1}^{|\Omega|} |b_{(i)}| \cdot |\delta_{(i)}| + |b_{(|\Omega|+1)}| \cdot \frac{|\Omega'|}{|\Omega|} \cdot \sum_{i=1}^{|\Omega|} |\delta_{(i)}| \\
&\leq \left(1 + \frac{|\Omega'|}{|\Omega|}\right) \sum_{i=1}^{|\Omega|} |b_{(i)}| \cdot |\delta_{(i)}| \\
&= \left(1 + \frac{|\Omega'|}{|\Omega|}\right) \langle b^*, \delta^*_{|\Omega|}\rangle,
\end{aligned} \tag{15}$$

where the vector $\delta^*_{|\Omega|}$ takes the absolute values of $\delta$ only at the positions where the vector $\delta$ has the $|\Omega|$ largest magnitudes and zeros elsewhere, that is,

$$\delta^*_{|\Omega|,i} = \begin{cases} |\delta_i|, & \text{if } |\delta_i| \geq |\delta_{(|\Omega|)}|, \\ 0, & \text{otherwise.} \end{cases}$$

For the vector $b^*$,

$$b_i^* = \begin{cases} |b_{(j)}|, & \text{if } \delta_i = \delta_{(j)} \text{ and } |b_{(j)}| \geq |b_{(|\Omega|)}|, \\ 0, & \text{otherwise.} \end{cases}$$

Combining (14) and (15), we have

$$\langle \Phi^T\varepsilon, \delta\rangle \leq \left(1 + \frac{|\Omega'|}{|\Omega|}\right) \langle b^*, \delta^*_{|\Omega|}\rangle. \tag{16}$$

Meanwhile, for any $\Omega^*$ we have

$$
\begin{aligned}
\|\Phi\delta\|_2^2 &\geq \|\Phi_{\Omega^*}\Phi_{\Omega^*}^+\Phi\delta\|_2^2 \\
&\geq \sigma_{\min,|\Omega^*|}^2 \cdot \|\Phi_{\Omega^*}^+\Phi\delta\|_2^2 \\
&= \sigma_{\min,|\Omega^*|}^2 \cdot \|\Phi_{\Omega^*}^+\Phi\delta_{\Omega^*} + \Phi_{\Omega^*}^+\Phi\delta_{\Omega^{*c}}\|_2^2,
\end{aligned}
\tag{17}
$$

where the set $\Omega^{*c}$ is the complement of the set $\Omega^*$, and the matrices $\Phi_{\Omega^*}$ and $\Phi_{\Omega^{*c}}$ are submatrices of the matrix $\Phi$ formed by taking columns whose indices are in $\Omega^*$ and $\Omega^{*c}$, respectively. As mentioned earlier, the matrix $\Phi_{\Omega^*}^+$ is a pseudo-inverse of $\Phi_{\Omega^*}$. The vector $\delta_{\Omega^*}$ (resp., $\delta_{\Omega^{*c}}$) takes only nonzero values when the index is in the set $\Omega^*$ (resp., $\Omega^{*c}$). In the above steps, the first inequality holds because the matrix $\Phi_{\Omega^*}\Phi_{\Omega^*}^+$ is a projection matrix. The second inequality is based on the definition of $\sigma_{\min,|\Omega^*|}^2$. The last step is a simple reformulation.

Note that in (17), $\Omega^*$ can be any subset of the indices. In the following, without loss of generality, we assume that the set $\Omega^*$ corresponds to the largest $|\Omega|$ magnitudes in the vector $\delta$, that is, $|\Omega^*| = |\Omega|$ and $\delta_{\Omega^*} = \delta_{|\Omega|}^*$. We then have

$$
\begin{aligned}
\|\Phi_\Omega^+\Phi\delta_\Omega + \Phi_{\Omega^*}^+\Phi\delta_{\Omega^{*c}}\|_2 &\geq \|\Phi_{\Omega^*}^+\Phi\delta_{\Omega^*}\|_2 + \|\Phi_{\Omega^*}^+\Phi\delta_{\Omega^{*c}}\|_2 \\
&\geq \|\delta_{\Omega^*}\|_2 - \sum_{k\in\Omega^{*c}} |\delta_k| \cdot \sup_{\ell\notin\Omega^*} \|\Phi_{\Omega^*}^+\phi_\ell\|_2 \\
&\geq \|\delta_{\Omega^*}\|_2 - \sum_{k=|\Omega|+1}^n |\delta_{(k)}| \cdot \sup_{\ell\notin\Omega^*} \|\Phi_{\Omega^*}^+\phi_\ell\|_2 \\
&\geq \|\delta_{\Omega^*}\|_2 - \frac{|\Omega'|}{|\Omega|} \cdot \|\delta_{|\Omega|}^*\|_1 \cdot \sup_{\ell\notin\Omega^*} \|\Phi_{\Omega^*}^+\phi_\ell\|_2 \\
&\geq \|\delta_{\Omega^*}\|_2 - \frac{|\Omega'|}{\sqrt{|\Omega|}} \cdot \|\delta_{|\Omega|}^*\|_2 \cdot \sup_{\ell\notin\Omega^*} \|\Phi_{\Omega^*}^+\phi_\ell\|_2 \\
&\geq \lambda(|\Omega|;|\Omega'|) \cdot \|\delta_{|\Omega|}^*\|_2.
\end{aligned}
\tag{18}
$$

In the above, the first and second steps are common manipulations. The third inequality takes $\Omega^*$ to be the subset of indices where $\delta_{\|\Omega\|}^*$ has nonzero entries. The fourth inequality is based on $\|\delta_{|\Omega|}^*\|_1/|\Omega| \geq \sum_{k=|\Omega|+1}^n |\delta_{(k)}|/|\Omega'|$. The fifth inequality is based on $\|\delta_{|\Omega|}^*\|_1 \leq \sqrt{|\Omega|} \cdot \|\delta_{|\Omega|}^*\|_2$. The last step recalls the definition of $\lambda(\cdot,\cdot)$. Combining (17) and (18), we have

$$
\|\Phi\delta\|_2^2 \geq \sigma_{\min,|\Omega|}^2 \cdot \lambda^2(|\Omega|;|\Omega'|) \cdot \|\delta_{|\Omega|}^*\|_2^2.
\tag{19}
$$

We now combine the above results and then maneuver back to the argument in the proof of Theorem 3.1. Combining (13), (16) and (19), we

have

$$\sigma_{\min,|\Omega|}^2 \cdot \lambda^2(|\Omega|;|\Omega'|) \cdot \|\delta_{|\Omega|}^*\|_2^2 \leq 2\left(1 + \frac{|\Omega'|}{|\Omega|}\right)\langle b^*, \delta_{|\Omega|}^*\rangle + \lambda_0(|\Omega| - |\Omega'|).$$

Let

$$\delta' = \frac{|\Omega|}{|\Omega| + |\Omega'|}\sigma_{\min,|\Omega|}^2 \cdot \lambda^2(|\Omega|;|\Omega'|) \cdot \delta_{|\Omega|}^*.$$

We have

$$\|\delta'\|_2^2 \leq 2\langle b^*, \delta'\rangle + \lambda_0 \cdot \frac{|\Omega|^2(|\Omega| - |\Omega'|)}{(|\Omega| + |\Omega'|)^2} \cdot \sigma_{\min,|\Omega|}^2 \cdot \lambda^2(|\Omega|;|\Omega'|).$$

The above is equivalent to

$$\|\delta' - b^*\|_2^2 \leq \|b^*\|_2^2 + \lambda_0 \cdot \frac{|\Omega|^2(|\Omega| - |\Omega'|)}{(|\Omega| + |\Omega'|)^2} \cdot \sigma_{\min,|\Omega|}^2 \cdot \lambda^2(|\Omega|;|\Omega'|),$$

which leads to

$$\|\delta'\|_\infty \leq \|b^*\|_\infty + \sqrt{c^2(b^*,|\Omega|) + \lambda_0 \cdot \frac{|\Omega|^2(|\Omega| - |\Omega'|)}{(|\Omega| + |\Omega'|)^2} \cdot \sigma_{\min,|\Omega|}^2 \cdot \lambda^2(|\Omega|;|\Omega'|)}.$$

Recalling the definitions of $\delta'$ and $b^*$, we have

$$
\begin{aligned}
\|\delta\|_\infty &\leq \bigg(c(\Phi^T\varepsilon, 1) \\
&\quad + \sqrt{c^2(\Phi^T\varepsilon,|\Omega|) + \lambda_0 \cdot \frac{|\Omega|^2(|\Omega| - |\Omega'|)}{(|\Omega| + |\Omega'|)^2} \cdot \sigma_{\min,|\Omega|}^2 \cdot \lambda^2(|\Omega|;|\Omega'|)}\bigg) \\
&\quad \times \left(\frac{|\Omega|}{|\Omega| + |\Omega'|}\sigma_{\min,|\Omega|}^2 \cdot \lambda^2(|\Omega|;|\Omega'|)\right)^{-1} \\
&\leq q'(|\Omega|;M).
\end{aligned}
$$
(20)

The above is equivalent to $\|x - x'\|_\infty \leq q'(|\Omega|;M)$. Using the same argument as in the last proof, we can argue that $\Omega \subset \Omega'$. Supposing $x_i \neq 0$, we have

$$|x'_i| \geq |x_i| - |x_i - x'_i| \geq |x_i| - q'(|\Omega|, M) > 0,$$

which implies that $\Omega \subset \Omega'$.

**Acknowledgment.** Comments from two anonymous referees led to significant reorganization of this paper. The original version [26] contains a lot more detail and is potentially more readable for newcomers.

SCHOOL OF INDUSTRIAL AND SYSTEMS ENGINEERING
GEORGIA INSTITUTE OF TECHNOLOGY
ATLANTA, GEORGIA 30332-0205
USA
E-MAIL: [xiaoming@isye.gatech.edu](xiaoming@isye.gatech.edu)
[xni@isye.gatech.edu](xni@isye.gatech.edu)